\documentclass[12pt,amsart]{article}
\usepackage{amsmath,amssymb,latexsym,epsfig,subfigure,multirow,setspace,mathrsfs}
\usepackage[utf8]{inputenc}
\usepackage[T1]{fontenc}
\newtheorem{definition}{Definition}[section]
\newtheorem{example}[definition]{Example}
\newtheorem{thm}[definition]{Theorem}
\newtheorem{lem}[definition]{Lemma}

\newtheorem{prop}[definition]{Proposition}

\newtheorem{remark}[definition]{Remark}

\newcommand\qed{\hfill\framebox[2.5mm]{}}

\newcommand{\Comment}[1]{}

\numberwithin{equation}{section}

\begin{document}

\title{Generalizations of wreath product identities via Garsia-Gessel bijections\footnote{This paper is dedicated to Bruce Sagan, for his profound insights on the topic of this paper, and his unconditional support to the author over the years.}}

\author{Tingyao Xiong\\ {\tt txiong@radford.edu}}
\date{}

\maketitle

{\begin{abstract} \textbf{Garsia-Gessel’s bijections have been proven to be very useful in obtaining multivariate  generating functions of permutation statistics. In 2011, Biaogioli and Zeng successfully derived four and six variate distributions via Garsia-Gessel’s bijections, by defining a ``Biaogioli-Zeng’s ordering” on the set of wreath product $\mathbb{Z}_r \wr S_n$. In this paper, we will prove that B-Z’s four-variate identities can be generalized to any positive-dominant ordering. We will also prove that B-Z’s six-variate distribution function can be significantly simplified under A-R’s ordering which was originally defined by Adin and Roichman in $2001$.}\\
\end{abstract}}


\section{Introduction}

Since Mc Mahon published his master volume {\it Combinatory Analysis} in $1915$, his esteemed successors Carlitz (\cite{Carlitz1}, \cite{Carlitz2}), Gessel $\&$  Garsia \cite{G-G}, and Stanley \cite{Stanley},  have contributed many important results about the generating functions and distribution identities with the following statistics or indices involved: \\
\begin{definition}\label{descendent}{\rm
Given a permutation $\pi=(\pi_1, \pi_2,\ldots,\pi_{n})\in S_n$, the descent set of $\pi$ is
\[
Des(\pi):=\{j\;|\;1\le j<n\;\text{ and }\; \pi_j>\pi_{j+1}\}
\]
then the {\em descent statistic} is 
\[
des(\pi):=\# Des(\pi)=\sum_{1\le j\le n-1}\chi(\pi_j>\pi_{j+1});
\]
the {\em major index} of $\pi$ is
\[
maj(\pi):=\sum_{j\in Des(\pi)}j=\sum_{1\le j\le n-1}j\cdot \chi(\pi_j>\pi_{j+1});
\]
and the {\em inversion number} of $\pi$ is 
\[
inv(\pi):=\#\{(i,j)\;|\;i<j\;\text{ and }\;\pi_i>\pi_j\}=\sum_{1\le i< j\le n}\chi(\pi_i>\pi_{j})
\]
}
\end{definition}
One famous identity of the two-variate distribution $(des(\pi), maj(\pi))$, that we will refer to at a later time, is called the {\em Carlitz Identity} (\cite{Carlitz2}):
\[
\frac{\sum_{\pi\in S_n}t^{des(\pi)}q^{maj(\pi)}}{ \prod_{j=0}^n (1-tq^j)}=\sum_{k\ge 0}[k+1]_q^n t^k
\]
where $[n]_q=1+q+q^2+\ldots+q^{n-1}$.\\

In 1979, Garsia and Gessel (\cite{G-G}) constructed two innovative bijections which paved a way of deriving the following two multivariate generating functions of joint distribution:
\begin{equation}\label{GG1}
\sum_{n\ge 0}\frac{\sum_{\pi\in S_n}p^{inv(\pi)}q^{maj(\pi)}t^{des(\pi)}}{\prod_{i=0}^n (1-tq^i)}\times \frac{u^n}{[n]_p!}=\sum_{k\ge 0}t^k \prod _{j=0}^ke[q^ju]_p
\end{equation}

and
\begin{equation}\label{GG2}
\frac{\sum_{\pi\in S_n}t_1^{des(\pi)} t_2^{des(\pi^{-1})}q_1^{maj(\pi)} q_2^{maj(\pi^{-1})}}{\prod_{i=0}^n (1-t_1q_1^i)(1-t_2q_2^i)}=\sum_{k_1\ge 0}\sum_{k_2\ge 0}t_1^{k_1} t_2^{k_2} \prod_{i\le k_1}\prod_{j\le k_2}\left.\frac{1}{1-uq_1^iq_2^j} \;\right\vert_{u^n},
\end{equation}
where $[n]_p!=[n]_p\times [n-1]_p\times \ldots  [1]_p$, and $e[u]_p=\sum_{n\ge 0}\frac{u^n}{[n]_p!}$.\\

Generalized definitions of colored permutation and the corresponding generation functions on the wreath product of $\mathbb{Z}_r$ by $S_n$ have been actively explored  in the past decades:
\begin{definition}\label{wreath}{\rm
For $r,n\in \mathbb{Z}_{>0}$, the wreath product of $\mathbb{Z}_r$ by $S_n$, or also  as known as $r-${\em colored permutation}, is defined as
\[
\mathbb{Z}_r \wr S_n=\{(\pi(1)^{c_1}, \pi(2)^{c_2},\ldots,\pi(n)^{c_n})\;|\;c_i\in [0, r-1], \pi\in S_n\}
\]}
\end{definition}
In 1993, Reiner (\cite{Reiner1},\cite{Reiner2}) generalized Garsia and Gessel's work to the set of {\it Coxter Group of type} $B_n$, which corresponds to the case $r=2$, or $2-$colored permutation. Reiner's research is often referred to by scholars as {\em signed permutation} since $\mathbb{Z}_2 \wr S_n$ could be represented as $\{(\pm\pi_1,\pm\pi_2,\ldots,\pm\pi_n)\}$.\\

When time moved into the post-millennium era, more and more analogous results on $\mathbb{Z}_r \wr S_n$ for any $r>0$ have been conquered. For instance,  different teams of scholars (\cite{Bagno},  \cite{B-B}, \cite{Chow&M}) have successfully extended Carlitz's identity to $\mathbb{Z}_r \wr S_n$ via various modifications on the traditional permutations statistics $maj(\pi)$ and $inv(\pi)$.\\

In 2011, following the proof framework  of Garsia and Gessel's bijection on the symmetric group $S_n$, Biagioli and Zeng (\cite{B&Z}) developed four-variate and six-variate joint distribution identities on $\mathbb{Z}_r \wr S_n$, which are akin to expressions (\ref{GG1}) and (\ref{GG2}). Biagioli and Zeng systematically updated the permutation statistics into the new forms to accommodate to the $r-$colored environment. In the following content, we will refer Biagioli and Zeng's work as B-Z's work, or B-Z's ordering.\\

The author would like to put some comments on B-Z's ordering. The table below lists the orders defined on $\mathbb{Z}_r \wr S_n$, that are cited accordingly from  \cite{A&R}, \cite{ST}, \cite{Reiner1} and \cite{B&Z}.
\begin{align*}
A-R:&\;\;1^{r-1}<\ldots<n^{r-1}<\ldots<1^1<\ldots<n^1<0<1<\ldots<n,& r>0; \\
ST:&\;\;0<1<\ldots<n<1^1<\ldots<n^1<\ldots<1^{r-1}<\ldots<n^{r-1},&r>0;\\
Re:&\;\;-n<-(n-1)<\ldots<-1<1<\ldots<n,&r=2;\\
B-Z:&\;\;n^{r-1}<\ldots<n^1<\ldots<1^{r-1}<\ldots<1^1<0<1<\ldots<n,&r>0.
\end{align*}
The order $A-R$ in the table above was originally studied by Adin and Roichman \cite{A&R}, which will be used in section $4$.\\

We can observe that the first three orders in the table are "color-based". In other words, these orders can be summarized as the following form:
\[
*^{r-1}\;\substack{<\\\text {(or }>)}\;*^{r-2}\ldots\substack{<\\\text{ (or }> )}\;*^1\;\substack{<\\\text{ (or }>)}\;*^0
\]
where $ *$ could be any number on $[n]$.

However, B-Z's ordering is "base-oriented", for its definition can be  abbreviated as
\[
n^\divideontimes<(n-1)^\divideontimes<\ldots<1^\divideontimes<0<1<2<\ldots<n,
\]
where $\divideontimes$ could be any positive color.\\

It is also noticeable that  Beck and Braun(\cite{Beck}), Davis and Segan (\cite{DS}) have generalized multiple identities, both on $S_n$ and $\mathbb{Z}_r \wr S_n$, through polyhedral geometrical methods.\\

In this paper, we will continue to apply Garsia-Gessel's bijections between the set of sequences of $r-$colored integers and the product set $\mathbb{Z}_r \wr S_n\times \mathscr{P}_n$.\\

Our paper is structured as follows: In section $2$, we will introduce definitions and notations to be used in the rest of this paper. The definitions in section 2 will be a similar but simpler version of the notations introduced in B-Z's paper. In section 3, we will prove that B-Z's four-variate joint distribution identity could be generalized to the cases under any positive-dominant ordering. In section 4, we will redefine a bipartite partition under the condition which is analogous to Garsia-Gessel's original construction. Under this new definition, we will apply A-R's ordering (the first order in the list above) to derive a new six-variate joint distribution which is significantly simpler than $(7.1)$ in \cite{B&Z}.

\section{More definitions and notations}
In previous section, we already introduced conventional notations such as $[n]_p!$ and $e[u]_p$. 

For $n\in \mathbb{N}$, define $[n]={1,2,\ldots,n}$. if both $n,m\in  \mathbb{Z}$ with $n\le m$, then let $[n,m]=\{n, n+1,\ldots,m\}$. Recall that $[n]_p!$ is already defined in section 1, $[0]_p!=1$. Now for $n=n_0+n_1+\ldots+n_k$ with $n_0, n_1,\ldots, n_k\ge 0$, the multinomial coefficient in $p$ form is
\[
\left[ \begin{tabular}{ cccc  }
\multicolumn{4}{c}{n} \\
$n_0$,& $n_1$,& \ldots, & $n_k$\\
\end{tabular}
\right]_p=\frac{[n]_p!}{[n_0]_p! \;[n_1]_p!\;\ldots \;[n_k]_p!}
\]

Given fixed $r,n \in \mathbb{N}$, we define $\mathfrak{G}(r,n)=\mathbb{Z}_r \wr S_n$, where the wreath product $\mathbb{Z}_r \wr S_n$ is as defined in Definition \ref{wreath}. 

\begin{definition}\label{col}{\rm
Given $r,n \in \mathbb{N}$, if $\gamma\in \mathfrak{G}(r,n)$, we use the notation as in (2.2) from \cite{B&Z}:
\begin{align*}
\gamma&=(\gamma(1),\gamma(2),\ldots,\gamma(n))=(\pi(1)^{c_1},\pi(2)^{c_2},\ldots,\pi(n)^{c_n})\\
&=\{<c_1,c_2\ldots,c_n>;\;\pi_\gamma=(\pi(1), \pi(2),\ldots,\pi(n))\}
\end{align*}
where $c_i\in[0, r-1]$. If $c_i=0$, then we will just write $\pi(i)^{c_i}=\pi(i)^0=\pi(i)$. Furthermore,
\begin{enumerate}
\item[(i)]the inverse of  $\gamma$ is $\gamma^{-1}=(\gamma^{-1}(1),\gamma^{-1}(2),\ldots,\gamma^{-1}(n)),\;\;\text{where}\;\;\gamma^{-1}(i)=\pi^{-1}(i)^{c_{\pi^{-1}(i)}},i\in[n]$.
\item[(ii)]The {\em color vector} of $\gamma$ is $Col(\gamma):=<c_1, c_2,\ldots,c_n>$,\;\;and the {\em color weight} of $\gamma$ is
\[
col(\gamma):=\sum_{i=1}^n c_i
\]
\end{enumerate}
}
\end{definition}

\begin{example}\label{inverse and color}{\rm
Suppose $\gamma =(5^2,\;2^3,\;4,\;3^1,\;7^2,\;1,\;6^1)\in \mathfrak{G}(4,7)$. So $\pi_\gamma=(5,2,4,3, 7,1,6)$, the color vector $Col(\gamma)=<2,3,0,1,2,0,1>\Rightarrow col(\gamma)=9$. We can calculate $\pi^{-1}_\gamma=(6,2,4,3,1,7,5)$, with the corresponding color distributions 

{\centering
$c_{\pi^{-1}(1)}=c_6=0,\;\; c_{\pi^{-1}(2)}=c_2=3,\;\;c_{\pi^{-1}(3)}=c_4=1,\;\;c_{\pi^{-1}(4)}=c_3=0,$\\
$c_{\pi^{-1}(5)}=c_1=2,\;\;c_{\pi^{-1}(6)}=c_7=1,\;\;c_{\pi^{-1}(7)}=c_5=2$.\\[1ex]
}
Thus $\gamma^{-1}=(6,\;2^3,\;4^1,\;3,\;1^2,\;7^1,\;5^2)$, with $Col(\gamma^{-1})=<0,3,1,0,2,1,2>$.
 }
\end{example}
\begin{remark}\label{symbolic}{\rm
In this paper,  $\mathfrak{G}(r,n)=\mathbb{Z}_r \wr S_n$ is a symbolic group. Although  we have defined $\gamma^{-1}$ on $\mathfrak{G}(r,n)$,,  there is no definition of the product operations among two elements in $\mathfrak{G}(r,n)$. That is, if $\alpha, \beta\in \mathfrak{G}(r,n)$, we have no definition on $\alpha\cdot \beta$.
}
\end{remark}
In section $1$, we have discussed several already-studied total orders, including B-Z's ordering, defined on the set of "colored" integers: 
\[
\mathfrak{C}(r,n):=\{n^{r-1},\;\dots,\;n^{1},\dots, 1^{r-1},\;\dots, 1^1, 0 ,\;1,\;\dots,\;n\}.
\]
The following definition gives a new classification of all the possible total orders on $\mathfrak{C}(r,n)$.
\begin{definition}\label{positive-do}{\rm
Given $r,n \in \mathbb{N}$, let $\mathfrak{C}(r,n)=\{j^c\;|\;j\in[0, n], c\in[0,r-1]\}$ be the set of colored entries of $\mathbb{Z}_r \wr S_n$. Recall that we have defined that for $j\in[0, n]$, $j^{c_j}=j$ if $c_j=0$. Then a total order defined on $\mathfrak{C}(r,n)$ is {\em positive-dominant} if
\begin{enumerate}
\item[(1)] $0<1<2<\ldots <n$;
\item[(2)] $i^0=i>j^c, \;\;\text{for any }\;i\in[0,n],\;j\in [n]\;\text{ and }\;c>0.$
\end{enumerate}
}
\end{definition}

\begin{example}\label{+dominant}{\rm
A-R's ordering and B-Z's ordering are positive-dominant. If $r=3$, $n=3$, then the following order is also positive-dominant:
\[
2^1<1^2<3^2<3^1<2^2<1^1<0<1<2<3
\]
}
\end{example}

\begin{remark}\label{positive dominant}{\rm
From the examples above we can see that under a positive-dominant ordering, the “positive” integers, or the numbers on $\mathfrak{C}(r,n)$ with zero color still follow the traditional ranking order; and these “positive” integers are greater than those integers with positive color values. For the author is concerned that altering the classical ordering $0<1<\ldots <n$ might cause unnecessary confusion in understanding the Garsia-Gassel’s construction which will be presented in Definitions \ref{colored permutation} and \ref{ball}. 
}
\end{remark}

Now we are ready to define the statistics of colored permutations, that are analogous to Definition \ref{descendent}.
\begin{definition}\label{colored per}{\rm
For a fixed pair $r,n \in \mathbb{N}$, given $\gamma=(\gamma(1),\gamma(2),\ldots,\gamma(n))\in \mathfrak{G}(r,n)$, let $\mathcal{O}$ be a total order defined on $\mathfrak{C}(r,n)$. If $\gamma(0)=0$, the descent set of $\gamma$ is
\[
Des_\mathcal{O}(\gamma):=\{j\;|\;j\in [0, n]\;\text{ and }\; \gamma(j)>\gamma(j+1)\};
\]
then the {\em descent statistic} of $\gamma$ is 
\[
des_\mathcal{O}(\gamma):=\# Des(\gamma)=\sum_{0\le j\le n-1}\chi(\gamma(j)>\gamma(j+1));
\]
the {\em major index} of $\gamma$ is
\[
maj_\mathcal{O}(\gamma):=\sum_{j\in Des(\gamma)}j=\sum_{0\le j\le n-1}j\cdot \chi(\gamma(j)>\gamma(j+1));
\]
and the {\em inversion number} of $\gamma$ is 
\[
inv_\mathcal{O}(\gamma):=\#\{(i,j)\;|\;1\le i<j\le n\;\text{ and }\;\gamma(i)>\gamma(j)\}=\sum_{1\le i< j\le n}\chi(\gamma(i)>\gamma(j));
\]
the {\em length} of $\gamma$ as defined in (\cite{Reiner1}, \cite{Bagno},  \cite{B&Z}) is
\[
len_\mathcal{O}(\gamma)=inv(\gamma)+\sum_{c> 0}(\pi(i)+c_i-1).
\]
In later sections, we sometimes omit the subscript $\mathcal{O}$, just use $des(\gamma)$, $maj(\gamma)$, etc, when there is no confusion.
}
\end{definition}
\begin{example}\label{example of inv}{\rm
Suppose $r=3$, $n=3$, we use the order $\mathcal{O}$ on $\mathfrak{C}(3,3)$ as mentioned in Example \ref{+dominant}:
\[
2^1<1^2<3^2<3^1<2^2<1^1<0<1<2<3
\] 
Then for $\gamma=(3^1,\; 1^2,\; 2^2)\in\mathfrak{G}(3,3)$, $Des(\gamma)=\{0,1\}\Rightarrow des(\gamma)=2$; $maj(\gamma)=0+1=1$, $inv(\gamma)=1$ for $\gamma$ only one inversion pair $(1,2)$, and finally
\[
len(\gamma)=inv(\gamma)+\sum_{c> 0}(\pi(i)+c_i-1)=1+(3+1-1)+(1+2-1)+(2+2-1)=9.
\]
}
\end{example}
The following two definitions are important for the following content.

\begin{definition}\label{inv on subset}{\rm
Suppose $\gamma=(\gamma(1),\gamma(2),\ldots,\gamma(n))\in \mathfrak{G}(r,n)$, and a total order $\mathcal{O}$ defined on $\mathfrak{C}(r,n)$. For a subset of colored entries $B=\{k_1^{c_{1}}, k_2^{c_{2}},\ldots, k_m^{c_{m}}\}\subseteq \mathfrak{C}(r,n)$, we define the inversion number of $\gamma$ on subset $B$ as
\[
Inv_{\mathcal{O}} (\gamma [B]):=\#\{(k_i,k_j)\;|\;1\le k_i<k_j\le n,\;\gamma(k_i),\gamma(k_j)\in B\;\text{ and }\;\gamma(k_i)>\gamma(k_j)\}.
\]
}
\end{definition}

\begin{example}\label{subset inv}{\rm
For $r=3$, $n=6$, let $\mathcal{Q}$ be the A-R's ordering on $\mathfrak{C}(3,6)$. Suppose
$B=\{1, 2^2, 3^1,4^2\}\subseteq \mathfrak{C}(3,6)$, and $\gamma =(1,\;6^1,\;3^1,\;4,\;2^2,\; 5^2)\in \mathfrak{G}(3,6)$, then
\[
inv_\mathcal{Q}(\gamma)=\#\{(1,2),(1,3),(1,5),(1,6),(2,3),(2,5),(2,6),(3,5),(3,6),(4,5),(4,6)\}=11; 
\]
while
\[
inv_\mathcal{Q}(\gamma[B])=\#\{(1,3),(1,5),(3,5)\}=3; 
\]
}
\end{example}
\begin{definition}\label{rank under an order}{\rm
Given a set $S$ with a total order $\mathcal{O}$, for any $a\in S$, we denote the order number of $a$ in $S$ under $\mathcal{O}$, from the smallest to the largest, as
\[
\mathcal{O}_ {S}(a)
\]
}
\end{definition}
\begin{example}\label{rank example}{\rm
Given a set $S=\{1,2,3\}$ with a total order $\mathcal{O}:1<2<3$, and another total order $\mathcal{P}:3<1<2$, we have
\[
\mathcal{O}_ {S}(2)=2,\text{ but }\;\;\mathcal{P}_ {S}(2)=3.
\]
}
\end{example}
\begin{definition}\label{reorder under another}{\rm
Given a colored permutation $\gamma=(\gamma(1),\gamma(2),\ldots,\gamma(n))\in \mathfrak{G}(r,n)$, choose a non-empty subset $B=\{\gamma(k_1), \gamma(k_2),\ldots,\gamma(k_m)\}\subseteq \mathfrak{C}(r,n)$. We construct a new colored permutation $\eta=(\eta(1),\eta(2),\ldots,\eta(n))$, such that
\begin{enumerate}
\item[i] if $\gamma(i)\in B$, define $\eta(i)$ such that $\mathcal{P}_B(\gamma(i))=\mathcal{O}_B(\eta(i))$;
\item[ii] if $\gamma(i)\not\in B$, define $\eta(i)=\gamma(i)$.
\end{enumerate}
Such a permutation is denoted as $\eta=\gamma[\mathcal{O}(B)]$. And obviously $inv_{\mathcal{O}}(\eta [B])= inv_{\mathcal{P}}(\gamma [B])$.
}
\end{definition}

\begin{example}\label{part change}{\rm
We will still use Example \ref{subset inv} to illustrate the concepts in Definition \ref{reorder under another}. For $r=3$, $n=6$, we denote A-R's ordering  as $Q$, and B-Z's ordering as
$P$ on set $\mathfrak{C}(3,6)$. Choose $\gamma =(1,\;6^1,\;3^1,\;4,\;2^2,\; 5^2)\in \mathfrak{G}(3,6)$, define $B=\{1, 2^2, 3^1\}$. Then under $ \mathcal{Q}$: $1>3^1>2^2$, but under $\mathcal{P}$: $1>2^2>3^1$. So we construct a new color permutation $\eta= \gamma[\mathcal{P}(B)]=(1,\;6^1,\;2^2,\;4,\;3^1,\; 5^2)$ by swapping the original positions of $3^1$ and $2^2$ for $\mathcal{Q}_B(3^1)= \mathcal{P}_B(2^2)$. And
\[
inv_\mathcal{Q}(\gamma[B])=\#\{(1,3^1), \{(1,2^2),( 3^1,2^2)\}=3, \;\;inv_\mathcal{P}(\eta[B])=\# \{(1,2^2),(1,3^1),( 2^2,3^1)\}=3
\]
So $inv_\mathcal{Q}(\gamma[B])= inv_\mathcal{P}(\eta[B])$.
}
\end{example}

\section{The four-variate distributions of $(des, maj, len, col)$}
In this section, we will extensively reference and reconstruct B-Z's results in sections 3, 4, and 5 of \cite{B&Z}. We will particularly reprove Lemma 4.3 in \cite{B&Z} (Lemma \ref{4.3} in this paper), which is the key step to prove Theorem \ref{four variate identity} (Theorem 5.1 in \cite{B&Z}), the most important result in this section. 

The sequence of colored integers defined as follows will play an important role in constructing the first bijection of this paper.
\begin{definition}\label{sequence of non-ne}{\rm
Given $r,n \in \mathbb{N}$, let 

\[
\mathbb{N}_0^{r,n}:=\{f=(f_1^{cf_1},\dots, f_n^{cf_n})\;|\;f_i\in \mathbb{N}\cup\{0\},\;cf_i\in [0,r-1],\; f_i^{cf_i}=f_i,\text{ if }cf_i=0\},
\]


Moreover, given $f\in \mathbb{N}_0^{r,n}$, define
\[
\max(f):=\max_{i\in [n]}(f_i),\;\;\text{ and } |f|:=\sum_{i=1}^n f_i
\]
}
\end{definition}
Note that in the definition above, we use $cf_i$'s to represent the colors of $f_i$'s, that are distinguished from the $c_i$'s in $\gamma(f)$ to be defined in Definitions \ref{colored permutation}.\\

In this section, unless otherwise specified, we always use $\mathcal{Q}$ to represent A-R's ordering, and $\mathcal{P}$ to represent B-Z's ordering.\\ 

The following definition is the canonical construction of colored permutation from $\mathbb{N}_0^{r,n}$ inspired by Garsis and Gessel's work (\cite{G-G}, page 291).
\begin{definition}\label{colored permutation}{\rm
Let $f\in\mathbb{N}_0^{r,n}$ and $0\le p_1< p_2<\ldots< p_k$ be the different values taken by $f_i$, $i\in [n]$. Please note that these $p_i$ values have no colors. Denote $A_{p_v}=\{i^{c_i}\;|\;f_i=p_v,v\in [n]\}$. Choose an order $\mathcal{O}$ on $\mathfrak{C}(r,n)$, list the entries in each $A_{p_v}$ in increasing order under $\mathcal{O}$, denoted as $\uparrow_\mathcal{O} A_{p_v}$. Then the colored permutation 
\[
\gamma_\mathcal{O}(f)=( \uparrow_\mathcal{O} A_{p_1}, \uparrow_\mathcal{O}  A_{p_2},\ldots, \uparrow_\mathcal{O}  A_{p_k})
\]
In future content, we may occasionally write $\gamma_\mathcal{O}(f)=\gamma(f)$ when the order is clear.
}
\end{definition}
\begin{example}\label{permutation under different color}{\rm
We will keep referring example \ref{subset inv}. For $r=3$, $n=6$, as mentioned before, let A-R's ordering be $\mathcal{Q}$, and B-Z's ordering be
$P$ on set $\mathfrak{C}(3,6)$. Now choose $f=(4^2, 3^1, 0, 2^2, 4^1,3^1)$. We rank the values of $f_i$ in increasing order as $0<2<3<4$, then we have $A_0=\{3\}$, $A_2=\{4^2\}$, $A_3=\{2^1, 6^1\}$, $A_4=\{1^2, 5^1\}$. Under different orders we have different colored permutations in $\mathfrak{G}(3,6)$:
\[
\gamma_\mathcal{Q}(f)=(3,4^2,2^1,6^1,1^2,5^1)=\{<0,2,1,1,2,1>;\pi=(3,4,2,6,1,5)\}
\] 
and
\[
\gamma_\mathcal{P}(f)=(3,4^2,6^1,2^1,5^1,1^2)=\{<0,2,1,1,1,2>;\sigma=(3,4,6,2,5,1)\}
\]
}
\end{example}
Besides Definition \ref{colored permutation}, we still need the following to reveal the connection between $\mathfrak{G}(r,n)$ and $\mathbb{N}_0^{r,n}$.
\begin{definition}\label{partition}{\rm
Given $n\in\mathbb{N}$, let $\mathscr{P}_n$ be the set of partitions of length $n$, that is
\[
\mathscr{P}_n =\{(\lambda_1,\lambda_2,\ldots,\lambda_n)\;|\; 0\le \lambda_1\le \lambda_2\le \ldots\le\lambda_n\text{ with all }\lambda_i\in \mathbb{N}\cup\{0\}\} 
\]
}
\end{definition}
\begin{example}\label{cp-partition}{\rm
Let $f\in\mathbb{N}_0^{r,n}$, let $\mathcal{O}$ be an order defined on $\mathfrak{C}(r,n)$, using the notation as in Definitions \ref{col} and \ref{colored permutation},
\[
\gamma_\mathcal{O}(f)= \{<c_1,c_2\ldots,c_n>;\;\pi\}.
\]
Then $(f_{\pi(1)},f_{\pi(2)},\ldots, f_{\pi(n)})$ is a partition.\\
For instance, in Example \ref{permutation under different color}, $f=(4^2, 3^1, 0, 2^2, 4^1,3^1)$ and 
\[
\gamma_\mathcal{Q}(f)=\{<0,2,1,1,2,1>;\pi=(3,4,2,6,1,5)\},
\]
then
\[
(f_{\pi(1)},f_{\pi(2)}, f_{\pi(3)},f_{\pi(4)}, f_{\pi(5)},f_{\pi(6)})=(0,2,3,3,4,4)
\]

Is a partition of length $6$.
}
\end{example}
The following proposition is a generalized version of Lemma 3.5 of \cite{B&Z}.
\begin{prop}\label{lambda}{\rm
Given $r,n \in \mathbb{N}$, $f\in\mathbb{N}_0^{r,n}$, under a positive-dominant ordering $\mathcal{O}$ on $\mathfrak{C}(r,n)$, construct $\gamma_\mathcal{O}(f)= \{<c_1,c_2\ldots,c_n>;\;\pi\}$ as in Definition \ref{colored permutation}, define sequence $\lambda(f):=(\lambda_1,\lambda_2,\ldots, \lambda_n)$ such that
\[
\lambda_i=f_{\pi(i)}-\#\{j\in Des\;|\;j\le i-1, \;\;i\in [n]\}
\]
Then $\lambda(f)$ is a partition.
}
\end{prop}
{\bf Proof.} Under any positive-dominant ordering $\mathcal{O}$ on $\mathfrak{C}(r,n)$, the proof of Proposition \ref{lambda} is identical to the proof of Lemma 3.5 of \cite{B&Z}. So interested readers can find the proof on page 542 \cite{B&Z}.\qed\\

The following definition creates important connections between the set of partitions and set of colored permutation.
\begin{definition}\label{partition and cor-p}{\rm
Given a partition $\lambda=(\lambda_1,\lambda_2,\ldots,\lambda_n)\in \mathscr{P}_n$, and a colored permutation $\gamma=\{<c_1,c_2\ldots,c_n>;\;\pi\}\in \mathfrak{G}(r,n)$, then define
\[
\lambda^\gamma=(\lambda_{\pi(1)}^{c_1}, \lambda_{\pi(2)}^{c_2},\ldots, \lambda_{\pi(n)}^{c_n})\in \mathbb{N}_0^{r,n}
\]
Moreover, we call a partition $\lambda$ is $\gamma-${\em compatible} if $\lambda_i<\lambda_{i+1}$ for all $i\in Des_\mathcal{O}(\gamma)$, where we define $\lambda_0=0$.
}
\end{definition}

Now we are ready to construct the bijection between $\mathbb{N}_0^{r,n}$ and $\mathfrak{G}(r,n)\times\mathscr{P}_n$. The following Lemma again, is a generalization of Lemma 3.8 \cite{B&Z} under any positive-dominant ordering.
\begin{lem}\label{1st bijection}{\rm
Given $r,n \in \mathbb{N}$, $f\in\mathbb{N}_0^{r,n}$, for a positive-dominant order $\mathcal{O}$ on $\mathfrak{C}(r,n)$, let $\gamma_\mathcal{O}(f)$ be as defined in Definition \ref{permutation under different color}, and  $\lambda(f)$ be as defined in Proposition \ref{lambda}, then the map
\[
\Phi: \mathbb{N}_0^{r,n}\mapsto \mathfrak{G}(r,n)\times\mathscr{P}_n \text { as } f\mapsto (\gamma_\mathcal{O}(f), \lambda(f))
\] 
is a bijection satisfying
\begin{enumerate}
\item[(i)] $\max(f)=\max(\lambda)+des_\mathcal{O}(\gamma (f))$
\item[(ii)]$|f|=|\lambda|+n\; des_\mathcal{O}(\gamma)-maj(\gamma)$.
\end{enumerate}
}
\end{lem}
{\bf Proof.} The proof of Lemma 3.8 \cite{B&Z} is correct for any positive-dominant ordering on $\mathfrak{C}(r,n)$. So we will not repeat the proof here. \qed

\begin{example}\label{bijection example 1}{\rm
Like in Example \ref{subset inv}, $f=(4^2, 3^1, 0, 2^2, 4^1,3^1) \in\mathbb{N}_0^{r,n}$, Then we have
\[
\Phi_\mathcal{Q} (f)=[(3,4^2,2^1,6^1,1^2,5^1), (0,1,2,2,2,2)]\in \mathfrak{G}(r,n)\times\mathscr{P}_n
\] 
and
\[
\Phi_\mathcal{P} (f)=[(3,4^2,6^1,2^1,5^1,1^2), (0,1,1,1,1,1)]\in \mathfrak{G}(r,n)\times\mathscr{P}_n
\]
Under $\mathcal{Q}$: $\max(f)=4$, $|\max(\lambda)|=2$, $des_\mathcal{Q}(\gamma(f))=\#\{(1,2),(4,5)\}=2$, so $\max(f)=\max(\lambda)+des_\mathcal{O}(\gamma (f))$;
$|f|=16$, $|\lambda|=9$, $n\; des_\mathcal{Q}(\gamma)=6\times 2=12$, $maj_\mathcal{Q}(\gamma)=1+4=5$. So $|f|=|\lambda|+n\; des_\mathcal{Q}(\gamma)-maj_{\mathcal{Q}}(\gamma)$.\\[1ex]
On the other hand, under $\mathcal{P}$: $des_\mathcal{P}(\gamma(f))=\#\{(1,2),(2,3)(4,5)\}=3$, $\max(f)=4$, $|\max(\lambda)|=1$,  so $\max(f)=\max(\lambda)+des_\mathcal{P}(\gamma (f))$.\\
Furthermore, $|f|=16$, $|\lambda|=5$, $n\; des_\mathcal{P}(\gamma)=6\times 3=18$, $maj_\mathcal{P}(\gamma)=1+2+4=7$. So $|f|=|\lambda|+n\; des_\mathcal{P}(\gamma)-maj_{\mathcal{P}}(\gamma)$.\\
}
\end{example}
\begin{remark}\label{remark about 1st bijection}{\rm
In Garsia and Gessel’s bijection, the image on $\mathscr{P}_n  $ was just 

$(f_{\pi(1)},f_{\pi(2)},\ldots, f_{\pi(n)})$. That is because when the permutations have no colors, the distribution of $des(\pi)$ is directly related to the values of $f_{\pi(i)}-f_{\pi(i+1)}$. However, on the wreath product $\mathfrak{G}(r,n)$, $des(\gamma)$ depends on the definition of ordering $\mathcal{O}$ on $\mathfrak{C}(r,n)$. That is why Biagioli and Zeng modified the partition sequence to $\lambda(f)$ as in Proposition \ref{lambda}.  At the same time, Example \ref{bijection example 1} also gives examples that the results of Lemma 3.5 and Proposition 3.8  of \cite{B&Z} can be generalized to Proposition \ref{lambda} and Lemma \ref{1st bijection} in this paper.
}
\end{remark}
Lemma \ref{1st bijection} immediately leads to the following result, which is also a generalizaztion to Proposition 3.12 \cite{B&Z}.

\begin{lem}\label{first identity}{\rm
Given $\sigma\in \mathfrak{G}(r,n)$, fix any positive-dominant ordering $\mathcal{O}$ on $\mathfrak{C}(r,n)$, we have
\[
\sum_{f\in \mathbb{N}_0^{r,n}\;|\gamma_\mathcal{O} (f)= \eta}t^{\max(f)}q^{\max(f)\cdot n-|f|}=\frac{t^{des_\mathcal{O}(\eta)}q^{maj_\mathcal{O}(\eta)}}{ \prod_{i=0}^{n-1}(1-tq^i)}
\]
}
\end{lem}
{\bf Proof.} This is a direct result from Lemma \ref{1st bijection}. \qed\\

In Garsia-Gessel’s bijection, there is an intermediate statistic $i(f)=\sum_{i<j}\chi(f_i<f_j)$. But now $f\in \mathbb{N}_0^{r,n}$, $i(f)$ will be updated to the version as follows.
\begin{definition}\label{inv(f) and col(f)}{\rm
Let $f=(f_1^{cf_1}, f_2^{cf_2},\ldots, f_n^{cf_n}) \in \mathbb{N}_0^{r,n}$, given an ordering $\mathcal{O}$ on $\mathfrak{C}(r,n)$, define
\[
Inv(f):=len_\mathcal{O}(\gamma(f)),\text{  and  }col(f)=\sum_{i=1}^ncf_i=col(\gamma(f))
\]
Let $\underline{n}=(n_0,n_1,\ldots,n_k)$ be a {\em composition} of $n$, that is, each $n_i\in \mathbb{N}\cup\{0\}$, and $n= n_0+n_1+\ldots+n_k$. Define
\[
\mathbb{N}_0^{r,n}(\underline{n}):=\{f\in \mathbb{N}_0^{r,n}\;|\;\#\{i:f_i=j\}=n_j\}
\]
}
\end{definition}
\begin{example}\label{composition 1}{\rm
Let $r=3$, $n=6$, let $\underline{6}=(3,2,1)$ be a composition of $6$.
Then both $f=(1^2,0,2^1, 1,0, 0)$ and $g=(0, 1^1,0,2^2, 1,0, 0)$ belong to $\mathbb{N}_0^{3,6}(\underline{6})$.
}
\end{example}
For a given composition $\underline{n}=(n_0,n_1,\ldots,n_k)$, we have the following proposition.
\begin{prop}\label{block form}{\rm
Let $r,n \in \mathbb{N}$, given $\underline{n}=(n_0,n_1,\ldots,n_k)$ a composition of $n$, for any positive-dominant order $\mathcal{O}$ on $\mathfrak{C}(r,n)$, define $\Gamma_{\underline{n}}\subseteq\mathfrak{G}(r,n)$ as
\[
\begin{array}{c c c c c}
\Gamma_{\underline{n}}=&\{( \underbrace{j_1,\dots,j_{n_0},} &\underbrace{j_{n_0+1}^{c_{n_0+1}},\dots,j_{n_0+n_1}^{c_{n_0+n_1}},}  &\ldots,& \underbrace{j_{n-n_k+1}^{c_{n-n_k+1}},\dots,j_{n_k}^{c_{n_k}}})\}\\
&\text{in }\uparrow_\mathcal{O}&\text{in }\uparrow_\mathcal{O}&\ldots&\text{in }\uparrow_\mathcal{O}
\end{array}
\]
Then $\Phi: \mathbb{N}_0^{r,n}(\underline{n})\mapsto  \Gamma_{\underline{n}}$ as $\Phi (f)=\gamma_\mathcal{O}(f)$ is a bijection, where $\gamma_\mathcal{O}(f)$ is as defined in Definition \ref{colored permutation}. 
}
\end{prop}
{\bf Proof.} To prove that $\Phi$ is a bijection, we only need to construct $\Phi^{-1}$. Given $\gamma\in \Gamma_{\underline{n}}$, then according to the definition of $\Gamma_{\underline{n}}$,
\[
\begin{array}{c c c c c}
\gamma=&\{( \underbrace{j_1,\dots,j_{n_0},} &\underbrace{j_{n_0+1}^{c_{n_0+1}},\dots,j_{n_0+n_1}^{c_{n_0+n_1}},}  &\ldots,& \underbrace{j_{n-n_k+1}^{c_{n-n_k+1}},\dots,j_{n_k}^{c_{n_k}}})\}\\
&\text{block 0 }&\text{block 1 }&\ldots&\text{block  }k
\end{array}
\]
Or equivalently,
\[
\begin{array}{c c c c c c}
\gamma=&\{( \underbrace{0,\dots,0} &\underbrace{c_{n_0+1},\dots,c_{n_0+n_1},}  &\ldots,& \underbrace{c_{n-n_k+1},\dots,c_{n_k}});&\pi\}\\
&\text{block 0 }&\text{block 1 }&\ldots&\text{block  }k&
\end{array}
\]
Now define $\Phi^{-1}$: for each entry  $\pi(i)^{c_i}$ in $\gamma$, suppose $\pi(i)^{c_i}$ locates in block $m$, $m\in [0, k]$, let the $\pi(i)_{th}$ entry of $f$ as $m^{c_i}$. Then for $i\in [0,k]$, such a $\Phi^{-1}(\gamma)$ has exactly $n_i$ many $i$ values with various colors. Thus $f=\Phi^{-1}(\gamma)\in\mathbb{N}_0^{r,n}(\underline{n})$.\qed
\begin{example}\label{want to die}{\rm
Let $r=3$, $n=6$ with a composition  $\underline{n}=(2,2,2)$. Under $\mathcal{Q}$, let $\gamma=(3,6,1^2,4^1,5^1,2)\in \Gamma_{\underline{n}}$. Block $0$ of $\gamma$ has two integers with color zero, then we define $f_3=f_6=0$; block $1$ has two colored integers $1^2$ and $4^1$, then put $1^2$ and $1^1$ as the first and fourth entries of $f$; Finally, $5^1$ and $2$ from block $2$ corresponds to $2^1$ on the fifth and $2$ on the second position of $f$. Combine the result above, we have
\[
f=(1^2, 2,0,1^1,2^1,0)\in \mathbb{N}_0^{r,n}(\underline{n}).
\]
}
\end{example}

The Lemma \ref{4.3}  below is an analogous result of Lemma 4.3 in \cite{B&Z}. However, we did not adopt the form in $(4.10)$ of \cite{B&Z} because we believe that $(\ref{Lemma 4.3})$ is a more intuitive formula to possibly develop a combinatorial interpretation of Lemma  \ref{4.3} in the future. We will prove that the conclusion in Lemma 4.3 in \cite{B&Z} can be extended to the cases of all positive-dominant orderings on $\mathfrak{C}(r,n)$.

\begin{lem}\label{4.3}{\rm
For $r,n \in \mathbb{N}$, given $\underline{n}=(n_0,n_1,\ldots,n_k)$ a composition of $n$, for any positive-dominant ordering $\mathcal{O}$ on $\mathfrak{C}(r,n)$, then
\begin{equation}\label{Lemma 4.3}
\sum_{f\in\mathbb{N}_0^{r,n}(\underline{n})}p^{inv_ \mathcal{O}(f)}a^{col(f)}=\prod_{i=n_0}^{n-1}\left(1+\sum_{j=1}^{r-1}a^jp^{j+i}\right) \left[ \begin{tabular}{ p{0.5cm}p{0.5cm}p{0.5cm}p{0.5cm}  }
\multicolumn{4}{c}{n} \\
$n_0$,&$n_1$,&\ldots,&$n_k$\\
\end{tabular}
\right]_p
\end{equation}
}
\end{lem}
{\bf Proof.} Let $\mathcal{P}$ be the B-Z's ordering. Lemma 4.3 in \cite{B&Z} already proved that expression $(\ref{Lemma 4.3})$ is true under $\mathcal{P}$. 

Now if $f\in\mathbb{N}_0^{r,n}(\underline{n})$, from Proposition \ref{block form} we have
\begin{equation}\label{blocks}
 \begin{array}{c c c c c}\gamma_\mathcal{P}(f)=& \underbrace{\{n_0 \text{ integers with no color }},&\underbrace{n_1 \text{  colored integers  }}, &\ldots,&\underbrace{n_k   \text{ colored integers }}\}\\
& \text{ positions of } 0&\text{ positions of } 1^c& \ldots&\text{ positions of } k^c\\
\end{array}
\end{equation}
Or equivalently,
\[
\gamma_\mathcal{P}(f)=\gamma(f)=\{(0,\dots,0,|c_{n_0+1},\dots,c_{n_0+n_1}|\;\ldots,\;|\;c_{n-n_k+1}, \dots,c_{n_k}), c_i\ge 0, i\in[n_0+1,n];\;\sigma \}
\]

In the rest of this proof, we will use $\gamma(f)$ to represent $\gamma_\mathcal{P}(f)$.

We note that expression $(\ref{blocks})$ is a general form of color-encoding a $f\in\mathbb{N}_0^{r,n}(\underline{n})$ for a specific $n-$ composition $\underline{n}=(n_0,n_1,\ldots,n_k)$. i.e. the entries in each block of $(\ref{blocks})$, only depend on the format of $f$, but  are independent of the choice of ordering on $\mathfrak{C}(r,n)$. Different ordering arranges the entries in each block differently. But the ranking difference occurred within a block will not affect the value of $inv(f)$. So the inversions of $f$ have to be obtained from pairs locating in different blocks. Specifically, we classify the inversions in $(\ref{blocks})$ into the following three cases. Without loss, in the list below, we always assume that the left number in the pair locates in a block left to the block of the right number: 
\begin{enumerate}
\item[(i)] $(b, s)$ with $b>s$, where both $b$ and $s$ are “positive”, or are integers having color zero.
\item[(ii)] $(b, s^c)$, where $b$ has color zero but $s$ has a positive color $c$.
\item[(iii)] $(b^{c_1}, s^{c_2})$, where  $c_1, c_2>0$ and $ b^{c_1}>_\mathcal{P} s^{c_2}$.
\end{enumerate}
So if we change $\mathcal{P}$ to any other positive-dominant ordering, only the inversions in case $(iii)$ as above are affected. This observation gives us the following construction:

Let $B=\{\sigma(i)^{c_i}\;|\;c_i>0\}$ be the set of  all colored entries in $(\ref{blocks})$. For any positive-dominant ordering $\mathcal{O}$, using the notations defined in Definition \ref{reorder under another}, construct $\delta=\gamma[\mathcal{O}(B)]$, then
\begin{enumerate}
\item[(1)] if $\gamma(i)\in B$, $\mathcal{P}_B(\gamma(i))=\mathcal{O}_B(\delta(i))$;
\item[(2)] if $\gamma(i)\not\in B$, $\gamma(i)=\delta(i)$.
\end{enumerate}
Then for any element $b\in B$, $\mathcal{P}_B(b)=\mathcal{O}_B(b)$, so $inv (\gamma[B])=inv(\delta[B])$. i.e. the inversions in case $(iii)$ stay the same on $\delta$. Changing the total ordering from $\mathcal{P}$ to $\mathcal{O}$ does not affect the inversions in cases $(i)$ and $(ii)\Rightarrow\;inv(\gamma(f))=inv(\delta)$. Note that this map $\gamma(f)\mapsto \delta$ is a bijection because for any $\delta\in \Gamma_{\underline{n}}$ arranged under $\mathcal{O}$, we can retreave $\gamma=\delta[\mathcal{P}(A)]$, where $A$ is the set of colored entries of $\delta$.
Furthermore, different ordering only shuffles the positions of colored entries without alternating the format of each colored entry, so $len_\mathcal{P}(\gamma(f))=len_\mathcal{O}(\delta)$. 

Now apply Proposition \ref{block form} to obtain
$g=\Phi^{-1}(\delta)\in\mathbb{N}_0^{r,n}(\underline{n})$, which gives a bijection $\Psi 
$ on $\mathbb{N}_0^{r,n}(\underline{n})$, such that if $f\in\mathbb{N}_0^{r,n}(\underline{n})$, $\Psi 
(f)=g$ satisfying 
\[
inv_\mathcal{P}(f)= inv_\mathcal{O}(g),\;\;\text{ and } col(f)= col(g)
\]
Since expression $(\ref{Lemma 4.3})$ is true under B-Z's ordering, the bijection $\Psi$ proves that $(\ref{Lemma 4.3})$ holds for any positive-dominant ordering on $\mathfrak{C}(r,n)$.\qed

\begin{example}\label{die tomorrow}{\rm
We will revisit Example \ref{want to die}: Let $r=3$, $n=6$ with a composition  $\underline{n}=(2,2,2)$. Under $\mathcal{Q}$, let $\delta=(3,6,1^2,4^1,5^1,2)\in \Gamma_{\underline{n}}$. Then $\Phi^{-1}(\delta)=f=(1^2, 2,0,1^1,2^1,0)\in \mathbb{N}_0^{r,n}(\underline{n})$ as shown in Example \ref{want to die}. $A=\{1^2,4^1,5^1\}$ is the set of colored entries in $\delta$. Under $\mathcal{Q}:\;1^2<4^1<5^1$; But under $\mathcal{P}:\;5^1<4^1<1^2$;
then construct $\gamma=\delta[\mathcal{P}(A)]= (3,6,5^1,4^1,1^2,2)$, then $\Phi^{-1}(\gamma)=g=(2^2, 2,0,1^1,1^1,0)$.

We can easily check that $\delta$ and $\gamma$ have exactly the same inversion pairs:
\begin{align*}
&inv_\mathcal{O}(\delta)=inv_\mathcal{P}(\gamma)=\#\{(1,3),(1,4),(1,5),(1,6), (2,3),(2,4),(2,5),(2,6)\}=8,\\
&col(\delta)=col(\gamma)=4, \;\;len_\mathcal{O}(\delta)=len_\mathcal{P}(\gamma)=19;
\end{align*}
}
\end{example}
\begin{definition}\label{general multinomial}{\rm
For $r,n \in \mathbb{N}$, define
\[
[n]_{r,a,p}=[\sum_{j=0}^{r-1}(a^jp^{j+n-1})^{\chi(j>0)}] \cdot [n]_p;\;\;\;\;[n]_{r,a,p}!=\prod_{k=1}^n \;[ k ]_{r,a,p}\;\;\;\;\text{and }\;\;\;e[u]_{r,a,p}=\sum_{n\ge 0}\frac{u^n}{[n]_{r,a,p}!}
\]
}
\end{definition}
Now we are ready to derive the four variate joint distribution of $(des, maj, len, col)$, which is analogous to expression (\ref{GG1}).
\begin{thm}\label{four variate identity}{\rm
For $r,n \in \mathbb{N}$, under any positive-dominant ordering $\mathcal{O}$ on $\mathfrak{C}(r,n)$, we have
\begin{equation}\label{4 variate identity}
\sum_{n\ge 0}\frac{\sum_{\gamma\in\mathfrak{G}(r,n)}t^{des_\mathcal{O}(\gamma)} q^{maj_\mathcal{O}(\gamma)} p^{len_\mathcal{O}(\gamma)} a^{col(\gamma)}
}{\prod_{i=0}^n(1-tq^i)} \frac{u^n}{[n]_{r,a,p}!}=\sum_{k\ge 0}t^k\prod_{j=0}^{k-1}e[q^ju]_p\cdot e[q^ku] _{r,a,p}
\end{equation}
}
\end{thm}
{\bf Proof.} Biagioli and Zeng have provided detailed proof of Theorem \ref{four variate identity} based on B-Z's ordering. B-Z’s proof was based on Proposition \ref{block form} and Lemma \ref{4.3}. which we have also successfully proved and extended to cases under any positive-dominant ordering. Thus the author would refer interested readers to find the proof in \cite{B&Z}, page 548.\qed

\section{The six-variate distributions}
In this section, we will use the method of bipartite partitions Garsia and Gessel invented in (\cite{G-G}), to derive a six-variate distribution of $(des(\gamma), des(\gamma^{-1}), maj(\gamma), maj(\gamma^{-1}), col(\gamma), col(\gamma)^{-1})$ which will be an analogous result of $(\ref{GG2})$. In Instead of using B-Z’s ordering, we will focus on A-R's ordering in this section. 

\begin{definition}\label{ball}{\rm
For $r,n \in \mathbb{N}$, let $\mathfrak{B}(r,n)$ be the set of colored bipartite partitions taking the forms   $\left[ \begin{array}{c} g \\ f \end{array}\right]= \left[ \begin{array}{ccc} g_1&\ldots&g_n \\ f_1^{cf_1}&\ldots&f_n^{cf_n} \end{array} \right]\in \mathscr{P}_n\times \mathbb{N}_0^{r,n}$, subjected to the condition that for each $i\in [n-1]$
\begin{enumerate}
\item[(a)] either $g_i<g_{i+1}$;
\item[(b)] or if $g_i=g_{i+1}$, $f_i^{cf_i}\le_{\mathcal{O}} f_{i+1}^{cf_{i+1}}$, under a given ordering $\mathcal{O}$ on the set of colored integers.
\end{enumerate}
}
\end{definition}
In this section, unless otherwise stated, we will always choose the ordering $\mathcal{O}$ to be the A-R's ordering: when  $i,j \in \mathbb{N}$,  
\begin{enumerate}
\item[(1)] if $c_i>c_j$, $i^{c_i}<j^{c_j}$; 
\item[(2)] if $c_i=c_j$, $i^{c_i}<(i+1)^{c_j}$.  
\end{enumerate}
And in the subsequent content, we will omit the subscript $\mathcal{O}$ from the indices like $des$, $maj$, etc.
\begin{example}\label{ball example}{\rm
Let $r=3$, $n=4$, the following two both belong to  $\mathfrak{B}(3,4)$.
\[
\left[\begin{array}{cccc} 0&1&1&1 \\ 2^{1}&2^{2}&2^{2}&3^1\end{array} \right],\;\;\; \left[ \begin{array}{cccc} 0&1&1&1 \\ 2^{2}&2^{1}&3^{1}&3\end{array} \right]
\]
}
\end{example}
The following Lemma reveals the relation between $\mathfrak{B}(r,n)$ and $\mathfrak{G}(r,n)$.
\begin{lem}\label{compatible}{\rm
Given $\left[ \begin{array}{c} g \\ f \end{array} \right]= \left[ \begin{array}{ccc} g_1&\ldots&g_n \\ f_1^{cf_1}&\ldots&f_n^{cf_n} \end{array} \right]\in\mathfrak{B}(r,n)$, let $\gamma(f)=\{<c_{1},\dots, c_{n}>;\;\pi \}$. Define $\mu=(f_{\pi(1)},\ldots, f_{\pi(n)})$ and $g=( g_1,\ldots,g_n)$. Then
\begin{enumerate}
\item[(1)] $\mu$ is $\gamma-$compatible;
\item[(2)]$g$ is $\gamma^{-1}-$compatible.
\end{enumerate}
}
\end{lem}
{\bf Proof.} Part $(1)$ is a well-known result due to the construction of  $\gamma(f)$ as in Definition \ref{colored permutation}. So we will only focus on proving part $(2)$: $g$ is $\gamma^{-1}-$compatible.

To prove that $g$ is $\gamma^{-1}-$compatible, we need to show that if $i\in Des(\gamma^{-1})$, $g_i<g_{i+1}$. Suppose
\[
\gamma^{-1}=\{<d_{1},\dots, d_{n}>;\;\pi^{-1}\}.
\]

If $0\in Des(\gamma^{-1})\rightarrow \gamma^{-1}(1)=k^d$, with $d>0$. That implies $\gamma(k)=1^d\rightarrow f_1>0$. So $g_1>0$ to meet the condition $(a)$ in Definition \ref{ball}.

If $0<i\in Des(\gamma^{-1})\rightarrow \pi^{-1}(i)^{d_i}>\pi^{-1}(i+1)^{d_{i+1}}$, there are the following possible cases:
\begin{itemize}
  \item[1.] $\pi^{-1}(i)<\pi^{-1}(i+1)$, $0\le d_i<d_{i+1}$, denote $\pi^{-1}(i)=u$ and $\pi^{-1}(i+1)=u+m$, $m>0$. Then $\gamma(u)=i^{d_i}$, $\gamma(u+m)=(i+1)^{d_{i+1}}$, so $f_i<f_{i+1}\& f_i^{d_i}>f_{i+1}^{d_{i+1}}$. Then $g_{i}<g_{i+1}$ by condition $(b)$ in Definition \ref{ball};
  \item[2.] $\pi^{-1}(i)>\pi^{-1}(i+1)$, denote $\pi^{-1}(i)=k+m$ and $\pi^{-1}(i+1)=k$, $m>0$.
\begin{itemize}
\item[\textbullet]if $d_i=d_{i+1}=0$, then $\pi(k)=i+1\;\&\;\pi(k+m)=i\Rightarrow f_{i}>f_{i+1}$ since both of them have color zero. Thus $g_{i}<g_{i+1}$ by condition $(b)$ in Definition \ref{ball};
\item[\textbullet] $d_i=0\&d_{i+1}>0$, then $\gamma(k)=(i+1)^{d_{i+1}}\;\&\;\gamma(k+m)=i\Rightarrow f_{i+1}\le f_i\Rightarrow f_i >f_{i+1}^{d_{i+1}} $, so $g_{i}<g_{i+1}$ by condition $(b)$ in Definition \ref{ball};
\item[\textbullet] $d_i=d_{i+1}=d>0$. Then $\gamma(k)=(i+1)^{d}\;\&\;\gamma(k+m)=i^ {d}\Rightarrow f_{i}>f_{i+1}\Rightarrow f_{i}^{d}>f_{i+1}^d$, then $g_{i}<g_{i+1}$ by condition $(b)$ in Definition \ref{ball};
\item[\textbullet] $d_i<d_{i+1}$. Then $\gamma(k)=(i+1)^{d_{i+1}}\;\&\;\gamma(k+m)=i^ {d_i}\Rightarrow f_{i}\ge f_{i+1}\Rightarrow f_{i}^{d_i}>f_{i+1}^{d_{i+1}}$, then $g_{i}<g_{i+1}$ by condition $(b)$ in Definition \ref{ball};\end{itemize}
\end{itemize}
Combine all $g$ is $\gamma^{-1}-$compatible.\qed\\
The Theorem below is an analogous result of Garsia-Gessel’s second bijection (Theorem 2.1, \cite{G-G}):
\begin{thm}\label{second bijection}{\rm
For $r,n \in \mathbb{N}$, there is a bijection between $\mathfrak{B}(r,n)$ and the triplets
\[
(\gamma,\;\lambda,\; \mu)
\]
Where $\gamma\in \mathfrak{G}(r,n)$, $\lambda, \mu\in \mathscr{P}_n$ satisfying that $\mu$ is $\gamma-$compatible and $\lambda$ is $\gamma^{-1}-$compatible.
}
\end{thm}

{\bf Proof. } Construct a map $\Phi: \mathfrak{B}(r,n)\mapsto (\gamma,\lambda, \mu)
$ as
\[
\Phi \left( \left[ \begin{array}{ccc} g_1&\ldots&g_n \\ f_1^{cf_1}&\ldots&f_n^{cf_n} \end{array} \right]\right)= (\gamma,\;\lambda,\; \mu)\
\]
Where $\lambda=(g_1, \ldots, g_n)$ and $\mu$ is as defined in Lemma \ref{compatible}. By Lemma \ref{compatible}, $\mu$ is $\gamma-$compatible and $\lambda$ is $\gamma^{-1}-$compatible. Then only need to prove that $\Phi$ is bijective.

Define $\Phi^{-1}: (\gamma,\;\lambda,\; \mu)\mapsto \mathfrak{B}(r,n)$ as
\[
g=\lambda,\;\;f=\mu^{\lambda^{-1}}
\]
Where $\mu^{\lambda^{-1}}$ is defined as in Definition \ref{partition and cor-p}.

It remains to show that $\Phi^{-1}(\gamma,\;\lambda,\; \mu)\in \mathfrak{B}(r,n)$.

Suppose $g_i=g_{i+1}$, since $g$ is $\gamma^{-1}-$compatible, $\gamma^{-1}(i)<\gamma^{-1}(i+1)$. Using the notations in Lemma \ref{compatible},  we check the following possible cases:
\begin{itemize}
\item[1.] $k=\pi^{-1}(i)< \pi^{-1}(i+1)=k+m$, $m>0$, and $d_i\ge d_{i+1}$. Then $\gamma(k)=i^{d_i}\&\gamma(k+m)=(i+1)^{d_{i+1}}$. So $f_i\le f_{i+1}\Rightarrow f_i^{d_i}\le f_{i+1}^{d_{i+1}}$ since $d_i\ge d_{i+1}$.
\item[2.] $k+m=\pi^{-1}(i)> \pi^{-1}(i+1)=k$, $m>0$, and $d_i> d_{i+1}$. Then $\gamma(k)=(i+1)^{d_{i+1}}\&\gamma(k+m)=i^{d_{i}}$. So $f_i\ge f_{i+1}\Rightarrow f_i^{d_i}< f_{i+1}^{d_{i+1}}$ since $d_i> d_{i+1}$.
\end{itemize}
Combine all the results above, $\Phi^{-1}(\gamma,\;\lambda,\; \mu)\in \mathfrak{B}(r,n)$ which implies that $\Phi$ is a bijection.\qed\\

Now we will give a six-variate distribution function of 
\[
(des(\gamma),\;des(\gamma^{-1},\;maj(\gamma), \;maj(\gamma^{-1}),\;col(\gamma), \;col(\gamma^{-1})),
\]
which is a result analogous to (\ref{GG2}).
\begin{thm}\label{six-variate}{\rm
For $r,n \in \mathbb{N}$, we have
\begin{align}\label{TX!}
&\frac{\sum_{\gamma\in\mathfrak{G}(r,n)}t_1^{des(\gamma)} t_2^{des(\gamma^{-1})}q_1^{maj(\gamma)} q_2^{maj(\gamma^{-1})} a^{col(\gamma)} b^{col(\gamma^{-1})}}{\prod_{i=0}^n (1-t_1q_1^i)(1-t_2q_2^i)}\\\nonumber
=&\sum_{k_1\ge 0}\sum_{k_2\ge 0}(t_1q_1^n)^{k_1} (t_2q_2^n)^{k_2} \prod_{i\le k_1}\prod_{j\le k_2}\prod_{m=0}^{r-1}\left.\frac{1}{1-uq_1^{-i}q_2^{-j}a^mb^m} \;\right\vert_{u^n},
\end{align}
}
\end{thm}
{\bf Proof.} In this proof, we denote $\mathfrak{B}(r,n)$ as $\mathfrak{B}$; $\mathscr{B}=\left[ \begin{array}{c} g \\ f \end{array}\right]$ and 
\[
 \mathfrak{B}(k_1, k_2)= \{\mathscr{B}\;|\;
\;\max(g)\le k_1,\; \max(f)\le k_2 \}
\]
We will firstly imitate the method Garsia and Gessel implemented in \cite{G-G}: for $r,n \in \mathbb{N}$, we can obtain $\mathscr{B}$ ``lexicographically” from the construction of $\mathfrak{B}$, and obtain the following identity which is analogous to $(2.7)$ in \cite{G-G}, by calculating the numbers of the colored columns $\left( \begin{array}{c} i \\ j^m \end{array}\right)$:  
\[
\sum_{\mathscr{B}\in \mathfrak{B}(k_1,k_2)} x^{|g|}y^{|f|}z^{col(f)}= \prod_{i\le k_1}\prod_{j\le k_2}\prod_{m=0}^{r-1}\left.\frac{1}{1-ux^{i}y^{j}z^m} \;\right\vert_{u^n}
\]

Now let $x=q_1^{-1}$, $y=q_2^{-1}$, $z=ab$ since $col(\gamma)=col(\gamma^{-1})$. Then the expression above becomes
\[
\sum_{\mathscr{B}\in \mathfrak{B}(k_1,k_2)} q_1^{-|g|}q_2^{-|f|}a^{col(f)} b^{col(f^{-1})}=
\prod_{i\le k_1}\prod_{j\le k_2}\prod_{m=0}^{r-1}\left.\frac{1}{1-uq_1^{-i}q_2^{-j}a^mb^m} \;\right\vert_{u^n}
\]
Multiplying by $(t_1q_1^n)^{k_1}(t_2q_2^n)^{k_2}$ and summing we obtain
\begin{align}\label{1st way}
&\sum_{\mathscr{B}\in \mathfrak{B}} \frac{t_1^{\max(g)}q_1^{n\cdot\max(g)-|g|}b^{col(f^{-1})} t_2^{\max(f)}q_1^{n\cdot\max(f)-|f|}a^{col(f)}}{(1-t_1q_1^n)(1-t_2q_2^n)}\\\nonumber
=&\sum_{k_1\ge 0}\sum_{k_2\ge 0}(t_1q_1^n)^{k_1}(t_2q_2^n)^{k_2}\prod_{i\le k_1}\prod_{j\le k_2}\prod_{m=0}^{r-1}\left.\frac{1}{1-uq_1^{-i}q_2^{-j}a^mb^m} \;\right\vert_{u^n}
\end{align}
On the other hand, by Theorem \ref{second bijection} and Lemma \ref{first identity}, we have
\begin{align}\label{BZ’s last}
&\;\;\sum_{\mathscr{B}\in \mathfrak{B}}t_1^{\max(g)}q_1^{n\cdot\max(g)-|g|}b^{col(f^{-1})} t_2^{\max(f)}q_1^{n\cdot\max(f)-|f|}a^{col(f)}\\\nonumber
=&\sum_{\gamma\in \mathfrak{G}(r,n)}\; \sum_{\substack{\lambda\in\mathscr{P}_n \\ \gamma^{-1}(g)=\lambda}}t_1^{\max(\lambda)}q_1^{n\;\max{\lambda}-|\lambda|}b^{col(\gamma^{-1})}\;\sum_{\substack{\mu\in\mathscr{P}_n \\ \gamma(f)=\mu}}t_2^{\max(\mu)}q_2^{n\;\max{\mu}-|\mu|}a^{col(\gamma)} \\\nonumber
=&\sum_{\gamma\in\mathfrak{G}(r,n)} \frac{t_1^{des(\gamma)} q_1^{maj(\gamma)} a^{col(\gamma)} t_2^{des(\gamma^{-1})} q_2^{maj(\gamma^{-1})} b^{col(\gamma^{-1})}}{\prod_{i\in [0,n-1]}(1-t_1q_1^i)(1-t_2q_2^i)},
\end{align}
since A-Z's ordering is a positive-dominant ordering.

Compare expressions $(\ref{1st way})$ and $(\ref{BZ’s last})$, we obtain $(\ref{TX!})$.\qed
\begin{remark}\label{final words}{\rm
Although both inherited from Garsia and Gessel’s second bijection, $(\ref{TX!})$ is significantly simpler than $(7.1)$ of \cite{B&Z}, for $\mathfrak{B}(r,n)$ in Definition \ref{ball} is constructed in a more straight forward way. To see the difference, consider a simple example: let $r=3, n=4$, in our construction, for the same $g=(1,\;1,\;1,\;1)\in\mathscr{P}_4$, we have only one $f=(2^2, 2^2, 2^1, 2^1)$ to make the colored biword $\left[ \begin{array}{c} g \\ f\end{array} \right]\in \mathfrak{B}(3,4)$; But in B-Z’s construction, all of the following colored bipartite partitions  have to be included:\\
\begin{tabular}{lll}
 $(1)\;\;f=(2^2, 2^2, 2^1, 2^1)$;\;\;&\;\; $(2)\;\;f=(2^2, 2^1, 2^1, 2^2)$;\;\;&$(3)\;\;f=(2^2, 2^1, 2^2, 2^1)$;\\
$(4)\;\;f=(2^1, 2^2, 2^1, 2^2)$;\;\;&\;\; $(5)\;\;f=(2^1, 2^1, 2^2, 2^2)$;\;\;&$(6)\;\;f=(2^1, 2^2, 2^2, 2^1)$.
\end{tabular}
}
\end{remark}
All the $f$'s listed above are very similar, the only difference is the distinct permutations on the color indices. These redundant repeats result in an unnecessarily complicated form in $(7.1)$ of \cite{B&Z}.
\section{Acknowledgements}
The author would like to express her extreme gratitude to Riccardo Biagioli for his endless patience and generous assistance that encouraged the author to strive forward.


\end{document}